\documentclass[12 pt]{article}
\usepackage[russian, english]{babel}
\usepackage[T1]{fontenc}
\usepackage{pifont,setspace,hyperref,amsthm,amsmath,amssymb,color,multirow,graphicx,subfigure, color }
\usepackage{enumitem}
\usepackage{inputenc}
\usepackage{amsmath,latexsym,amsthm}
\usepackage{graphicx}
\usepackage{amssymb}

\begin{document}

\begin{center}
\textbf{\Large \textsc {A boundary-value problem for a mixed type equation involving hyper-Bessel fractional differential operator and  Hilfer's bi-ordinal fractional derivative}}\\[0.2 cm]

\bigskip
\textbf{Karimov E.~T.~\footnote{V.I.Romanovskiy Institute of Mathematics, Tashkent, Uzbekistan. E-mail: erkinjon@gmail.com},
	Ruzhansky~M.~\footnote{Department of Mathematics: Analysis, Logic and Discrete Mathematics, Ghent University, Ghent, Belgium and School of Mathematical Sciences,  Queen Mary University of London, United Kingdom. E-mail:  michael.ruzhansky@ugent.be},
	Toshtemirov B.~H.~\footnote{V.I.Romanovskiy Institute of Mathematics, Tashkent, Uzbekistan. E-mail: toshtemirovbh@gmail.com}.}

\end{center}

\textbf{Abstract.} { In a rectangular domain, a boundary-value problem is considered for a mixed-type equation with a regularized Caputo-like counterpart of hyper-Bessel differential operator and the bi-ordinal Hilfer's fractional derivative. Using the method of separation of variables, Laplace transform, a unique solvability of the considered problem has been established. Moreover, we have found explicit solution of initial problems for differential equations with the bi-ordinal Hilfer's derivative and regularized Caputo-like counterpart of hyper-Bessel differential operator with the non-zero starting point.}

\medskip
\noindent \textbf{MSC 2010: }35M10, 34A12

\smallskip
\noindent \textbf{Keywords:} Sub-diffusion equation; fractional wave equation; bi-ordinal Hilfer's derivative; hyper-Bessel fractional differential operator; boundary-value problems; Laplace transform.

\section{Introduction}

 The study of fractional order differential equations has been attracting many scientists because of its adequate and interesting applications in modeling of real-life problems related to several fields of science  \cite{L1}-\cite{L5}. Initial-value problems (IVPs) and boundary-value problems (BVPs) involving the Riemann-Liouville and Caputo derivatives attract  most interest (see, for instance, \cite{L6},  \cite{L7}, \cite{L8}).  Especially, studying IVPs and BVPs for the sub-diffusion, fractional wave equations are well-studied (see \cite{sad}, \cite{ruzh}, \cite{karruzh}). BVPs for mixed type equations are also an interesting target for many authors (see \cite{L9}-\cite{L13}).

 Introducing a generalized Riemann-Liouville fractional derivatives (it is called  Hilfer's deriva\-tive) has opened a new gate in the research of fractional calculus (\cite{L14}-\cite{L16}). Therefore, one can find several works devoted to studying this operator in various problems \cite{L17}, \cite{L18}.  We also note that in 1968,  M.~M.~Dzhrbashyan and A.~B.~Nersesyan introduced the following integral-differential operator \cite{Ldn}
\begin{equation}\label{DN}
D_{0x}^{\sigma_n}g(x)=I_{0x}^{1-\gamma_n}D_{0x}^{\gamma_{n-1}}...D_{0x}^{\gamma_{1}}D_{0x}^{\gamma_{0}}g(x), ~ ~ n\in\mathbb{N}, ~ x>0,
\end{equation}
which is more general than Hilfer's operator. Here $I_{0x}^{\alpha}$ and $D_{0x}^{\alpha}$ are the Riemann-Liouville fractional
integral and the Riemann-Liouville fractional derivative of order $\alpha$ respectively (see Definition 2.1),  $\sigma_n\in(0, n]$ which is defined by
$$
\sigma_n=\sum\limits_{j=0}^{n}\gamma_j-1>0, ~ \gamma_j\in(0, 1].
$$
There are some works \cite{bag}, \cite{bag2}, related with this operator. New wave of researches involving this operator might appear due to the translation of original work \cite{Ldn} in FCAA \cite{Lfca}.
 
 In addition, from announcing the concept of hyper-Bessel fractional differential derivative by I. Dimovski \cite{L21}, several articles have been published dedicated to studying problems containing this type of operators (see \cite{L22}-\cite{L26}).   For instance, fractional diffusion equation and wave equation were widely investigated in different domains in \cite{L19}-\cite{L20}.

In this work, we investigate a boundary value problem for a mixed equation involving the sub-diffusion equation with Caputo-like counterpart of a hyper-Bessel fractional differential operator and the fractional wave equation with Hilfer's bi-ordinal derivative in a  rectangular domain. The theorem about the uniqueness and existence of the solution is proved.

{The rest of the paper is organized as follows: In Preliminaries section we provide necessary information on Mittag-Leffler functions (Section 2.1), hyper-Bessel functions (Section 2.2.), bi-ordinal Hilfer's fractional derivatives (Section 2.3) and on differential equation involving bi-ordinal Hilfer's fractional derivatives (Section 2.4). Auxiliary result is formulated in Theorem 2.2.  In Section 3, we formulate  the main problem and state our main result in Theorem 3.1.  In Appendix one can find detailed arguments of the proof of Theorem 2.1.}

\section{Preliminaries}

In this section we present some definitions and auxiliary  results related to generalized Hilfer's derivative and fractional hyper-Bessel differential operator which will be used in the sequel. We start recalling the definition of the Mittag-Leffler function.

\subsection{Important properties of the Mittag-Leffler function}

The two parameter Mittag-Leffler (M-L) function is an entire function given by
\begin{equation}\label{e1}
E_{\alpha, \beta}(z)=\sum_{k=0}^{\infty}\frac{z^k}{\Gamma(\alpha{k}+\beta)}, \, \, \, \alpha>0, \, \beta\in\mathbb{R}.
\end{equation}

\textbf{Lemma 2.1} (see \cite{L3}) Let $\alpha<2, \, \beta\in\mathbb{R}$ and $\frac{\pi\alpha}{2}<\mu<min\{\pi, \pi\alpha\}$. Then the following estimate holds
$$
|E_{\alpha, \beta}(z)|\leq\frac{M}{1+|z|}, \, \, \, \mu\leq|argz|\leq\pi, \, \, |z|\geq0.
$$
Here and in the rest of the paper, $M$ denotes a positive constant.

In \cite{L27}, the following estimate for $E_{\delta, 1}(-t^\delta)$ in the the form  $g_{\delta}(t)\leq{E}_{\delta, 1}(-t^\delta)\leq{f_{\delta}(t)}$ was given without proof, where
$$
g_{\delta}(t)=\frac{1}{t^\delta\Gamma(1-\delta)+1}\sim\frac{t^{-\delta}}{\Gamma(1-\delta)}\sim{E}_{\delta, 1}(-t^\delta), \, t \to \infty,
$$
$$
f_{\delta}(t)=\frac{1}{1+\frac{t^\delta}{\Gamma(1+\delta)}}\sim(1-\frac{t^{\delta}}{\Gamma(1+\delta)}\sim{E}_{\delta, 1}(-t^\delta), \, t \to +0.
$$

Recently,  the lower and upper bounds of the Kilbas-Saigo function, which is a  more general form of $E_{\delta, 1}(-t^\delta)$, were announced. Also  the following preposition  about the bounds for two parameter M-L function  is given in \cite{L28}.

\textbf{Proposition 2.1} (\cite{L28}) For every $\alpha\in(0,1]$, \, $\beta>\alpha$ and $x\geq0$ one has
$$
\frac{1}{\Big(1+\sqrt{\frac{\Gamma(1-\alpha)}{\Gamma(1+\alpha)}}x\Big)^2}\leq{E_{\alpha, \alpha}(-x)}\leq\frac{1}{\Big(1+\frac{\Gamma(1+\alpha)}{\Gamma(1+2\alpha)}x\Big)^2}
$$
and
$$
\frac{1}{1+\frac{\Gamma(\beta-\alpha)}{\Gamma(\beta)}x}\leq\Gamma(\beta){E_{\alpha, \beta}(-x)}\leq\frac{1}{1+\frac{\Gamma(\beta)}{\Gamma(\beta+\alpha)}x}.
$$

The Laplace transform of M-L function is given in the following lemma.

\textbf{Lemma 2.2.} (\cite{L8}) For any $\alpha>0, \, \, \beta>0$  and $\lambda\in\mathbb{C}$, we have
$$
\mathcal{L}\{t^{\beta-1}E_{\alpha, \beta}(\lambda{t}^\alpha)\}=\frac{s^{\alpha-\beta}}{s^\alpha-\lambda}, \, \, (Re(s)>\mid{\lambda}\mid^{1/\alpha}),
$$
where the Laplace transform of  a function $f(t)$ is defined by

$$
\mathcal{L}\{f\}(s):=\int_{0}^{\infty}e^{-st}f(t)dt.
$$

\textbf{Lemma 2.3.}\label{lemma} If $\alpha\leq{0}$  and $\beta\in\mathbb{C}$, then the following recurrence formula holds:
$$
E_{\alpha, \beta}(z)=\frac{1}{\Gamma(\beta)}+zE_{\alpha, \alpha+\beta}(z).
$$
This lemma was proved by R. K. Saxena in 2002 \cite{L29}.

Later, we use the properties of a Wright-type function studied by A. Pskhu \cite{L30}, defined as
$$
e^{\mu, \delta}_{\alpha, \beta}(z)=\sum_{n=0}^{\infty}\frac{z^n}{\Gamma(\alpha{n}+\mu)\Gamma(\delta-\beta{n})}, \, \, \, \alpha>0, \, \, \alpha>\beta.
$$

M-L function can be determined by Wright-type function as a special case $E_{\alpha, \beta}(z)=e^{\beta, 1}_{\alpha, 0}(z)$.  So, we can record some properties of M-L function which can be reduced from the Wright-type function's properties.

\textbf{Lemma 2.4.} (\cite{L30}) If $\pi\geq|argz|>\frac{\pi\alpha}{2}+\varepsilon, \, \, \, \varepsilon>0$,  then the following relations are valid for $z \to \infty$:
$$
\mathop {\lim }\limits_{\mid{z}\mid \to  \infty} E_{\alpha, \beta}(z)=0,
$$
$$
\mathop {\lim }\limits_{\mid{z}\mid \to  \infty} zE_{\alpha, \beta}(z)=-\frac{1}{\Gamma(\beta-\alpha)}.
$$

\subsection{Regularized Caputo-like counterpart of the hyper-Bessel fractional differential operator}

\textbf{Definition 2.1.} (\cite{L8}) The Riemann-Liouville fractional integral $I^{\alpha}_{a+}f(t)$ and derivative  $D^{\alpha}_{a+}f(t)$ of order $\alpha$ are defined by
$$
I^{\alpha}_{a+}f(t)=\frac{1}{\Gamma(\alpha)}\int_{a}^{t}\frac{f(\tau)d\tau}{(t-\tau)^{1-\alpha}},
$$

$$
D^{\alpha}_{a+}f(t)=\left(\frac{d}{dt}\right)^nI^{n-\alpha}_{a+}f(t), \, \, \, \,   n-1\leq\alpha<n,
$$
where $\Gamma(\alpha)$ is Euler's gamma-function.

\textbf{Definition 2.2}. The Erdelyi-Kober (E-K) fractional integral of a function  $f(t)\in{C_\mu}$ with arbitrary parameters $\delta>0, \gamma\in{\mathbb{R}}$ and $\beta>0$ is defined as (\cite{L8})
$$
I^{\gamma, \delta}_{\beta; a+}f(t)=\frac{t^{-\beta(\gamma+\delta)}}{\Gamma(\delta)}\int_{a}^{t}(t^\beta-\tau^\beta)^{\delta-1}\tau^{\beta\gamma}f(\tau)d(\tau^\beta),
$$
which can be reduced up to a weight to  $I^{q}_{a+}f(t)$ (Riemann-Liouville fractional integral) at $\gamma=0$ and $\beta=1$, and Erdelyi-Kober fractional derivative of $f(t)\in{C}_{\mu}^{(n)}$ for $n-1<\delta\leq{n}, n\in\mathbb{N}$, is defined by
$$
D_{\beta, a+}^{\gamma, \delta}f(t)=\prod_{j=1}^{n}\left(\gamma+j+\frac{t}{\beta}\frac{d}{dt}\right)\big(I_{\beta, a+}^{\gamma+\delta, n-\delta} f(t)\big),
$$
where $C_{\mu}^{(n)}$  is the  weighted space of continuous functions defined as
$$
C_{\mu}^{(n)}=\left\{f(t)=t^p \tilde{f(t)}; \, \, \tilde{f} \in{C}^{(n)}[0, \infty)\right\}, \, \, \, C_{\mu}=C_{\mu}^{(0)} \, \,  \textrm{with} \, \, \mu\in\mathbb{R}.
$$

\textbf{Definition 2.3.} Regularized Caputo-like counterpart of the hyper-Bessel fractional differential operator for  $\theta<1$,  $0<\alpha\leq1$ and $t>a\geq{0}$ is defined in terms of the E-K fractional order operator

\begin{equation}\label{e2}
^{C}\Big(t^{\theta}\frac{d}{dt}\Big)^\alpha{f}(t)=(1-\theta)^\alpha{t}^{-\alpha(1-\theta)}D_{1-\theta, a+}^{-\alpha, \alpha}\left(f(t)-f(a)\right)
\end{equation}
or in terms of the hyper-Bessel differential (R-L type) operator

\begin{equation}\label{eq3}
^{C}\left(t^{\theta}\frac{d}{dt}\right)^\alpha{f}(t)=\Big(t^{\theta}\frac{d}{dt}\Big)^\alpha{f}(t)-\frac{f(a)\Big(t^{(1-\theta)}-a^{(1-\theta)}\Big)^{-\alpha}}{(1-\theta)^{-\alpha}\Gamma(1-\alpha)},
\end{equation}
where
$$
\left(t^\theta\frac{d}{dt}\right)^\alpha{f}(t)=\left\{ \begin{gathered}
(1-\theta)^\alpha t^{-(1-\theta)\alpha}I_{1-\theta, a+}^{0, -\alpha} f(t) \, \, \,  \textrm{if}  \, \,  \,  \theta<1, \hfill \cr
(\theta-1)^\alpha t^{-(1-\theta)\alpha}I_{1-\theta, a+}^{-1, -\alpha} f(t)  \, \, \,  \textrm{if} \, \,  \,  \theta>1,  \hfill \cr
\end{gathered}  \right.
$$
is a hyper-Bessel fractional differential operator (\cite{L25}).

From (\ref{eq3}) for $a=0$ we obtain the definition presented in (\cite{L25}) and also Caputo FDO is the particular case of Caputo-like counterpart hyper-Bessel operator at $\theta=0$.

\textbf{Theorem 2.1.}  Assume that the following conditions hold:

$\bullet$ $\tau\in{C}[0,1]$ such that $\tau(0)=\tau(1)=0$ and $\tau'\in{L}^2(0,1)$,

$\bullet$ $f(\cdot, t)\in{C}^3[0,1]$ and $f(x, \cdot)\in{C}_\mu[a, T]$ such that

 $f(0,t)=f(\pi,t)=f_{xx}(0,t)=f_{xx}(1,t)=0$, and $\dfrac{\partial^4}{\partial{x}^4}f(\cdot,t)\in{L}^1(0,1).$

Then, in  $\Omega=\{0<x<1, \, \, a<t<T \}$,  the problem of finding the solution of the equation
$$
^C\Big(t^\theta\frac{\partial}{\partial{t}}\Big)^{\alpha}u(x,t)-u_{xx}(x,t)=f(x,t),
$$
satisfying the conditions
$$
u(0,t)=0, u(1,t)=0, \, \, \, a\leq{t}\leq{T},
$$
$$
u(x,0)=\tau(x), \, \, \, \, 0\leq{x}\leq{1},
$$
has a unique solution given by
\begin{equation}\label{e4}
u(x,t)=\sum_{k=1}^{\infty}\Big[\tau_kE_{\alpha, 1}\Big(-\frac{(k\pi)^2}{p^\alpha}(t^{p}-a^{p})^\alpha\Big)+G_k(t)\Big]\sin(k\pi{x}),
\end{equation}
where $p=1-\theta$ and
$$
\begin{array}{l}
\displaystyle{G_k(t)=\frac{1}{p^\alpha\Gamma(\alpha)}\int_{a}^{t}\big(t^p-\tau^p\big)^{\alpha-1}f_k(\tau)d(\tau^p)}\\
\\
\displaystyle{-\frac{(k\pi)^2}{p^{2\alpha}}\int_{a}^{t}\big(t^p-\tau^p\big)^{2\alpha-1}E_{\alpha, 2\alpha}\Big[-\frac{(k\pi)^2}{p^\alpha}(t^p-\tau^p)^{\alpha}\Big]f_k(\tau)d(\tau^p),}
\end{array}
$$
$$
\tau_k=2\int_{0}^{1}\tau(x)\sin(k\pi{x})dx, \, \, \, \, 
f_k(t)=2\int_{0}^{1}f(x,t)\sin(k\pi{x})dx, \, \, \, \, k=1,2,3,...
$$

In fact, for $a=0$  Theorem 2.1  implies the result of (\cite{L22}) (see Theorem 3.1). For the detailed proof of Theorem 2.1 see appendix.
\subsection{ Bi-ordinal Hilfer's fractional derivative}
\textbf{Definition 2.4.} Hilfer's derivative $D^{\alpha, \mu}_{a+}$ of order $\alpha \, \, \,  (n-1<\alpha\leq{n},   \, \, \, n\in\mathbb{N})$  of type $\mu \, \,  (0\leq\mu\leq1)$ is defined by \cite{L14}:
\begin{equation}\label{e5}
D^{\alpha, \mu}_{t}f(t)=I^{\mu(n-\alpha)}_{0+}\left(\frac{d}{dt}\right)^nI^{(1-\mu)(n-\alpha)}_{0+}f(t).
\end{equation}
Then in \cite{L20}, V. M. Bulavatsky considered generalized Hilfer's derivative in the form
$$
D^{(\alpha, \beta) \mu}_{t}f(t)=I^{\mu(1-\alpha)}_{0+}\frac{d}{dt}I^{(1-\mu)(1-\beta)}_{0+}f(t),
$$
here  $0<\alpha, \beta\leq{1}$,  $0\leq\mu\leq{1}$.

In the same way, one can present Hilfer's bi-ordinal fractional derivative of orders \\ $\alpha \, \,  (n-1<\alpha\leq{n})$, $\beta \, \, (n-1<\beta\leq{n})$ and of type $\mu\in[0,1]$  by the following relation:

\begin{equation}\label{e6}
D^{(\alpha, \beta) \mu}_{t}f(t)=I^{\mu(n-\alpha)}_{0+}\left(\frac{d}{dt}\right)^nI^{(1-\mu)(n-\beta)}_{0+}f(t).
\end{equation}

In general, (\ref{e6}) is also preserved as (\ref{e5}) in terms of its interpolation concept. Specifically, when $\mu=0$, (\ref{e6}) gives the Riemann-Liouville fractional derivative of $\beta$ order and for $\mu=1$, the bi-ordinal Hilfer's derivative expresses the Caputo fractional derivative of order $\alpha$.

Similarly, generalized Hilfer's  fractional derivative of orders $\gamma\in(1,2]$, and $\beta\in(1,2]$ and type $\mu\in[0,1]$ can be written as a special case of (\ref{e6}) for $n=2$:
\begin{equation}\label{e7}
D^{(\gamma, \beta)\mu}_tf(t)=I^{\mu(2-\gamma)}_{0+}\left(\frac{d}{dt}\right)^2I^{(1-\mu)(2-\beta)}_{0+}f(t).
\end{equation}

Here we present the formula for the Laplace transform of (\ref{e7}) which will be used later:
$$
  \mathcal{L}\{D^{(\alpha, \beta)\mu}_{t}f(t)\}=s^{\beta+\mu(\alpha-\beta)}\mathcal{L}\{f(t)\}-
 $$
  \begin{equation}\label{e8}
  -s^{1-\mu(2-\alpha)}\left[I_{0+}^{(1-\mu)(2-\beta)}f(t)|_{t \to 0+}\right]-s^{-\mu(2-\alpha)}\left[\frac{d}{dt}I_{0+}^{(1-\mu)(2-\beta)}f(t)|_{t \to 0+}\right].
\end{equation}

\subsection{Differential equation involving bi-ordinal Hilfer's derivative}

Let us consider the following problem:

\textit{Find a solution of the equation}
\begin{equation}\label{e9}
D^{(\gamma, \beta)\mu}_{t}u(t)+\lambda{u}(t)=f(t),  \, \, \,   (1<\gamma, \beta\leq2, \, \, 0\leq\mu\leq{1}),
\end{equation}
\textit{satisfying the initial conditions}
\begin{equation}\label{e10}
\mathop {\lim }\limits_{t \to  0+}{I^{(1-\mu)(2-\beta)}_{0+}}u(t)=\xi_0,
\end{equation}
\begin{equation}\label{e11}
\mathop {\lim }\limits_{t \to  0+}\frac{d}{dt}{I^{(1-\mu)(2-\beta)}_{0+}}u(t)=\xi_1,
\end{equation}
\textit{where} $f(t)$ \textit{is a given function} $f\in{L^1}(0, \infty)$, and $\lambda, \, \, \xi_0, \, \,  \xi_1=const.$

\textbf{Theorem 2.2.} If  $f\in{C^1}(0, +\infty)$, then  the problem (\ref{e9})-(\ref{e11}) has an unique solution represented by
$$
u(t)=\xi_0t^{(\beta-2)(1-\mu)}E_{\delta, \delta+\mu(2-\gamma)-1}(-\lambda{t}^\delta)+\xi_1t^{\mu+(\beta-1)(1-\mu)}E_{\delta, \delta+\mu(2-\gamma)}(-\lambda{t}^\delta)+
$$
\begin{equation}\label{e12}
	+\int_{0}^{t}(t-\tau)^{\delta-1}E_{\delta, \delta}(-\lambda(t-\tau)^{\delta})f(\tau)d\tau,
\end{equation}
where $\delta=\beta+\mu(\gamma-\beta)$.

{\sf Proof:}
 Applying the Laplace transform (\ref{e9}) by means of (\ref{e8}) and considering initial conditions (\ref{e10}), (\ref{e11}) yield
\begin{equation}\label{e13}
\mathcal{L}\{u\}=\frac{\xi_0s^{1-\mu(2-\gamma)}+\xi_1s^{-\mu(2-\gamma)}+\mathcal{L}\{f\}}{s^{\beta+\mu(\gamma-\beta)}+\lambda},
\end{equation}
where $\mathcal{L}\{u\}$ and $\mathcal{L}\{f\}$ are the Laplace transform of functions $u$ and $f$, respectively.

According to Lemma 2.2, the Laplace transform of the Mittag-Leffler function \cite{L3}, \cite{L8} as follows
$$
\mathcal{L}^{-1}\{\frac{s^{1-\mu(2-\gamma)}}{s^{\beta+\mu(\gamma-\beta)}+\lambda}\}=t^{\beta-2+\mu(2-\beta)}E_{\beta+\mu(\gamma-\beta), \beta-1+\mu(2-\beta)}(-\lambda{t}^{\beta+\mu(\gamma-\beta)}),
$$
$$
\mathcal{L}^{-1}\{\frac{s^{-\mu(2-\gamma)}}{s^{\beta+\mu(\gamma-\beta)}+\lambda}\}=t^{\beta-1+\mu(2-\beta)}E_{\beta+\mu(\gamma-\beta), \beta+\mu(2-\beta)}(-\lambda{t}^{\beta+\mu(\gamma-\beta)}),
$$

$$
\mathcal{L}^{-1}\{\frac{\mathcal{L}\{f\}}{s^{\beta+\mu(\gamma-\beta)}+\lambda}\}=\int_{0}^{t}(t-\tau)^{\beta-1+\mu(\gamma-\beta)}E_{\beta+\mu(\gamma-\beta), \beta+\mu(\gamma-\beta)}(-\lambda(t-\tau)^{\beta+\mu(\gamma-\beta)})f(\tau)d\tau,
$$
where $\mathcal{L}^{-1}$ is an inverse Laplace transform operator.

Considering above evaluations and after applying the inverse Laplace transform to (\ref{e13}), we can write the solution of (\ref{e9})-  (\ref{e11}) in the form  (\ref{e12}). Theorem 2.2 is proved.

\section{Formulation of a problem and main result}

Let us consider the following equation
\begin{equation}\label{e14}
f(x, t) = \left\{ \begin{gathered}
	^C\left( {{t^\theta }\frac{\partial }{{\partial t}}} \right)^\alpha u\left( {x,t} \right) -
	{u_{xx}}\left( {x,t} \right),\,\,\,\,(x,t)\in{\Omega_1},  \hfill \cr
	{D^{(\gamma, \beta)\mu}_{t}u}\left( {x,t} \right) - {u_{xx}}\left( {x,t} \right),\,\,\,\,(x,t)\in{\Omega_2},  \hfill \cr
\end{gathered}  \right.
\end{equation}
in a domain $\Omega=\Omega_1\cup\Omega_2\cup{Q}$. Here
$\Omega_1=\{(x,t): 0<x<1,\, \,  \, a<t<b\}$,  \, \, \, \,  \\   $\Omega_2=\{(x,t): 0<x<1, 0<t<a\}$,
$Q=\{(x,t): 0<x<1,   t=a \}$, \, $a, b\in\mathbb{R^+}$ such that  $a<b$, $0<\alpha\leq1$, \,  $\theta<1$, \, $1<\gamma, \beta<2$, \, $0\leq\mu\leq1$, $f(x,t)$ is a given function, $^C\Big(t^\theta\frac{\partial}{\partial{t}}\Big)^\alpha$ is the regularized Caputo-like counterpart of the hyper-Bessel operator defined as in (\ref{e2}), $D^{(\gamma, \beta)\mu}_{t} $ is the bi-ordinal Hilfer's derivative defined as in (\ref{e7}).

\textbf{Problem}. Find a solution of (\ref{e14}) in $\Omega$, satisfying regularity conditions
$$
u(x,t)\in{C}(\Omega)\cap{C^1}(\Omega)\cap{C^2}(\Omega_2), \, \, \, ^C\Big(t^\theta\frac{\partial}{\partial{t}}\Big)^\alpha{u}(x,t)\in{C}(\Omega_1), \,
$$
$$
u_{xx}(x,t)\in{C^2}(\Omega_1), \, \, \, \, \, D_{t}^{(\gamma, \beta)\mu}u(x,t)\in{C(\Omega_2)},
$$
and the boundary-initial conditions
\begin{equation}\label{e15}
u(0,t)=0, \, \, \, \, \,  \,   0\leq{t}\leq{b},
\end{equation}
\begin{equation}\label{e16}
u(1,t)=0,  \, \,  \,   \, \, \, 0\leq{t}\leq{b},
\end{equation}
\begin{equation}\label{e17}
\mathop {\lim }\limits_{t \to  0+}{I^{(1-\mu)(2-\beta)}_{0+}}u(x,t)=\varphi(x), \, \, \, \, 0\leq{x}\leq1,
\end{equation}
as well as the gluing conditions

\begin{equation}\label{e18}
\mathop {\lim }\limits_{t \to  a-}{I^{(1-\mu)(2-\beta)}_{0+}}u(x,t)=\mathop {\lim }\limits_{t \to  a+}u(x,t), \, \, \, \, \,  0\leq{x}\leq{1},
\end{equation}

\begin{equation}\label{e19}
\mathop {\lim }\limits_{t \to  a-}\frac{d}{dt}{I^{(1-\mu)(2-\beta)}_{0+}}u(x,t)=\mathop {\lim }\limits_{t \to  a+}(t-a)^{1-(1-\theta)\alpha}u_t(x,t), \, \, \, \, \,  0<x<1,
\end{equation}
here  $\varphi(x)$ is a given function.

The key motivation to formulate this problem is a possible application in diffusion-wave processes, which will be described by the mixed type equation as Eq. (\ref{e14}) \cite{KDSc}. Moreover, fractional derivatives used in the mixed equation have more general and also specific character. Therefore, we separately studied related Cauchy problems for ordinary differential equations of fractional order.

The intention of this paper is to prove the uniqueness and existence of the solution to the problem (\ref{e14})-(\ref{e19}), as stated in the following theorem.

\textbf{Theorem 3.1.} If the following conditions

1)  $0<\beta+\mu(\gamma-\beta)\leq1$,   \\

2)  $\varphi\in{C}[0,1]$ and $\varphi'\in{L}^2[0,1],$ \\

3)  $f(\cdot,t)\in{C}^3[0,1]$ and $f(x,\cdot)\in{C}_\mu$ such that $f(0,t)=f(1,t)=f_{xx}(0,t)=f_{xx}(1,t)=0$ and
$\dfrac{\partial^4}{\partial{x}^4}f(\cdot, t)\in{L^1(0,1)}$
hold, then there exists a unique solution of the considered problem (\ref{e14})-(\ref{e19}).

{\sf Proof:}

First we introduce new notations:
\begin{equation}\label{e20}
\mathop {\lim }\limits_{t \to  0+}\frac{d}{dt}{I^{(1-\mu)(2-\beta)}_{0+}}u(x,t)=\psi(x), \, \, \,  0<x<1,
\end{equation}

\begin{equation}\label{e21}
	\mathop {\lim }\limits_{t \to  a+}u(x,t)=\tau(x),  \, \, \, \, 0\leq{x}\leq1,
\end{equation}
here $\tau(x)$ and $\psi(x)$ are unknown functions to be found later.

Using the method of separation of variables for solving the homogeneous equation corresponding to (\ref{e14}), i.e. searching solution as $u(x,t)=T(t)X(x)$ and considering  (\ref{e15}) and (\ref{e16}) in homogeneous case, yield the following problem:

\begin{equation}\label{e22}
	X''(x)+\lambda{X}(x)=0, \, \, \ \ X(0)=0, \, \, X(1)=0.
\end{equation}
It is obvious that (\ref{e22}) is a Sturm-Liouville problem on finding eigenvalues and eigenfunctions and it has the following solution:

\begin{equation}\label{e23}
\lambda_k=(k\pi)^2, \, \, \, X_k(x)=\sin(k\pi{x}),  \, \, \, k=1, 2, 3, ....
\end{equation}

Using the fact that the system of eigenfunctions $\{X_k\}$ in (\ref{e23}) forms an orthogonal basis in $L^2(0,1)$ \cite{L31}, we look for the solution $u(x,t)$ and given function $f(x,t)$ in the form of series expansions as follows:

\begin{equation}\label{e24}
u(x,t)=\sum_{k=1}^{\infty}u_k(t)\sin(k\pi{x}),	
\end{equation}
\begin{equation}\label{e25}
	f(x,t)=\sum_{k=1}^{\infty}f_k(t)\sin(k\pi{x}),	
\end{equation}
where $u_k(t)$ is unknown function to be found, $f_k(t)$ is known and given by
$$
f_k(t)=2\int_{0}^{1}f(x,t)\sin(k\pi{x})dx.
$$

Substituting (\ref{e24}) and (\ref{e25}) into equation (\ref{e14}) in $\Omega_1$ and considering initial condition (\ref{e20}) gives the following fractional differential equation
$$
^C\left(t^\theta\frac{d}{dt}\right)^{\alpha}u_k(t)+(k\pi)^2u_k(t)=f_k(t)
$$
with initial condition
$$
u_k(a+)=\tau_k,
$$
where $\tau_k$ is the coefficient of series expansion of $\tau(x)$ in terms of orthogonal basis (\ref{e22}), i.e.,
$$
\tau_k=2\int_{0}^{1}\tau(x)\sin(k\pi{x})dx.
$$
After finding the solution of this problem, then considering (\ref{e23}) we can write the solution of (\ref{e14}) in  $\Omega_1$ satisfying the conditions (\ref{e15}), (\ref{e16}) and  (\ref{e21}) stated in (\ref{e4}).

Now by using the solution (\ref{e4}), we evaluate $(t-a)^{1-(1-\theta)\alpha}u_t(x,t)$:
$$
(t-a)^{1-(1-\theta)\alpha}u_t(x,t)=\sum_{k=1}^{\infty}\left[-\frac{(k\pi)^2}{p^{\alpha-1}}\tau_k \, E_{\alpha, \alpha}\Big(-\frac{(k\pi)^2}{p^\alpha}(t^p-a^p)^{\alpha}\Big)+(t-a)^{1-p\alpha}G_k(t)\right]\sin(k\pi{x}).
$$

Considering above given evaluations we obtain the following relation on $Q$ deduced from $\Omega_1$ as $t\to{a+}$:

\begin{equation}\label{e26}
\mathop {\lim }\limits_{t \to  a+}(t-a)^{1-(1-\theta)\alpha}u_t(x,t)=\sum_{k=1}^{\infty}\left[-\frac{(k\pi)^2}{\Gamma(\alpha)p^{\alpha-1}}\tau_k \right]\sin(k\pi{x}).
\end{equation}

Now we establish another relation on $Q$ which will be reduced from $\Omega_2$.

According to variable separation method, considering (\ref{e24}), (\ref{e25}) and initial conditions (\ref{e17}), (\ref{e20}), we obtain the following problem finding a solution of equation
$$
D_{t}^{(\gamma, \beta)\mu}W(t)+\lambda_k{W}(t)=f(t)
$$
and satisfying the initial conditions
$$
\mathop {\lim }\limits_{t \to  0+}{I^{(1-\mu)(2-\beta)}_{0+}}W(t)=\varphi_k,
$$
$$
\mathop {\lim }\limits_{t \to  0+}\frac{d}{dt}{I^{(1-\mu)(2-\beta)}_{0+}}W(t)=\psi_k.
$$

It is obvious that (\ref{e12}) is the solution for above-given problem. Hence, using the solution (\ref{e12}) and taking  into account (\ref{e24})  we write the solution of (\ref{e14}) in $\Omega_2$ satisfying (\ref{e17}) and (\ref{e20}) as
\begin{equation}\label{e27}
u(x,t)=\sum_{k=1}^{\infty}W_k(t)\sin(k\pi{x}),
\end{equation}
where
$$
W_k(t)=\varphi_kt^{(\beta-2)(1-\mu)}E_{\delta, \delta+\mu(2-\gamma)-1}(-\lambda_k{t}^\delta)+\psi_kt^{\mu+(\beta-1)(1-\mu)}E_{\delta, \delta+\mu(2-\gamma)}(-\lambda_k{t}^\delta)+
$$
$$
+\int_{0}^{t}(t-\tau)^{\delta-1}E_{\delta, \delta}\left[-\lambda_k(t-\tau)^\delta\right]f_k(\tau)d\tau,
$$
here $\delta=\beta+\mu(\gamma-\beta)$ and $\varphi_k$ is coefficient of the series expansion of  $\varphi(x)$, i.e.,
$$
\varphi_k=2\int_{0}^{1}\varphi(x)\sin(k\pi{x})dx
$$
and $\psi_k$ is not known yet.

Now using (\ref{e27}) we simplify $\mathop {\lim }\limits_{t \to  a-}{I^{(1-\mu)(2-\beta)}_{0+}}W_k(t)$ and $\mathop {\lim }\limits_{t \to  a-}\dfrac{d}{dt}{I^{(1-\mu)(2-\beta)}_{0+}}W_k(t)$ in the following formulas

\begin{equation}\label{e28}
\begin{gathered}
\mathop {\lim }\limits_{t \to  a-}{I^{(1-\mu)(2-\beta)}_{0+}}W_k(t)=\varphi_k \, E_{\delta, 1}(-\lambda_k{a}^\delta)+\psi_k a\, E_{\delta, 2}(-\lambda_k{a}^\delta)\\
+\int_{0}^{a}(a-s)^{\delta+q-1}E_{\delta, \delta+q}\big[-\lambda_k(a-s)^\delta\big]f_{k}(s)ds,
\end{gathered}
\end{equation}
\begin{equation}\label{e29}
\begin{gathered}
\mathop {\lim }\limits_{t \to  a-}\frac{d}{dt}{I^{(1-\mu)(2-\beta)}_{0+}}W_k(t)=-\varphi_k\lambda_k{a}^{\delta-1} \, E_{\delta, \delta}(-\lambda_k{a}^\delta)+\psi_k \, E_{\delta, 1}(-\lambda_k{a}^\delta)\\
+\int_{0}^{a}(a-s)^{\delta+q-2}E_{\delta, \delta+q-1}\big[-\lambda_k(a-s)^\delta\big]f_{k}(s)ds.
\end{gathered}
\end{equation}
After substituting (\ref{e28}) and (\ref{e21}) into gluing condition (\ref{e18})  and substituting (\ref{e29}), (\ref{e26}) into the gluing condition (\ref{e19}), we obtain the following the system of linear algebraic equations with respect to $\tau_k$ and $\psi_k$:

\begin{equation}\label{e30}
 \left\{ \begin{gathered}
  \varphi_k \, E_{\delta, 1}(-\lambda_k{a}^\delta)+\psi_k a\, E_{\delta, 2}(-\lambda_k{a}^\delta)
		+\int_{0}^{a}(a-s)^{\delta+q-1}E_{\delta, \delta+q}\big[-\lambda_k(a-s)^\delta\big]f_{k}(s)ds=\tau_k,\\
\varphi_k\lambda_k{a}^{\delta-1} \, E_{\delta, \delta}(-\lambda_k{a}^\delta)-\psi_k \, E_{\delta, 1}(-\lambda_k{a}^\delta)-\\
\int_{0}^{a}(a-s)^{\delta+q-2}E_{\delta, \delta+q-1}\big[-\lambda_k(a-s)^\delta\big]f_{k}(s)ds=\frac{\lambda_k}{\Gamma(\alpha)p^{\alpha-1}}\tau_k.\\
 \end{gathered}  \right.
\end{equation}

From (\ref{e30}), we find $\psi_k$ and $\tau_k$:
\begin{equation}\label{e31}
\psi_k=\frac{B}{\Delta}\varphi_k+\frac{C}{\Delta},
\end{equation}
$$
\tau_k=\Big(E_{\delta, 1}(-\lambda_k{a}^\delta)+\frac{B}{\Delta}E_{\delta, 2}(-\lambda_k{a}^\delta)\Big)\varphi_k +\frac{C}{\Delta}E_{\delta, 2}(-\lambda_k{a}^\delta)+
$$
\begin{equation}\label{e32}
+\int_{0}^{a}(a-s)^{\delta+q-1}E_{\delta, \delta+q}\big[-\lambda_k(a-s)^\delta\big]f_{k}(s)ds,
\end{equation}
here
$$
\Delta=E_{\delta, 1}(-\lambda_k{a}^\delta)+
\frac{\lambda_k{a}}{\Gamma(\alpha){p}^{\alpha-1}}E_{\delta, 2}(-\lambda_k{a}^\delta),
$$

$$B=\frac{-\lambda_k{p}^{1-\alpha}}{\Gamma(\alpha)}E_{\delta, 1}(-\lambda_k{a}^\delta)+\lambda_k{a}^{\delta-1}E_{\delta, \delta}(-\lambda_k{a}^\delta),$$
$$
C=\frac{-\lambda_k{p}^{1-\alpha}}{\Gamma(\alpha)}\int_{0}^{a}(a-s)^{\delta+q-1}E_{\delta, \delta+q}\big[-\lambda_k(a-s)^\delta\big]f_{k}(s)ds-\int_{0}^{a}(a-s)^{\delta+q-2}E_{\delta, \delta+q-1}\big[-\lambda_k(a-s)^\delta\big]f_{k}(s),
$$
here $\lambda_k=(k\pi)^2$.

The system of linear equations (\ref{e30}) is equivalent to the considered problem in terms of existing the solution. For that reason, if  $\Delta\neq{0}$ for any $k$,  (\ref{e30}) has only one solution or the considered problem's solution is unique if it exists. Therefore, we show that $\Delta$ is not equal to zero for any $k$.

By using Lemma 2.4,  the behavior of $\Delta$ at $k \to \infty$ can be written as:
$$
 \mathop {\lim }\limits_{k \to  \infty}\Delta=\mathop {\lim }\limits_{\mid{z}\mid \to  \infty}\big[E_{\delta, 1}(z)+\frac{p^{1-\alpha}}{\Gamma(\alpha)a^{\delta-1}}zE_{\delta, 2}(z)\big]=\frac{p^{1-\alpha}}{\Gamma(\alpha)\Gamma(2-\delta)a^{\delta-1}}
$$

This proves that $\Delta\neq{0}$ for sufficiently large $k$.

According to Theorem 3.1 and using Proposition 2.1  we can show that
$$\Delta=E_{\delta,1}(-\lambda_k{a}^\delta)+\frac{\lambda_k{a}}{\Gamma(a)p^{\alpha-1}}E_{\delta, 2}(-\lambda_k{a}^\delta)\geq\frac{1}{1+\Gamma(1-\delta)\lambda_k{a}^\delta}+\frac{1}{1+\Gamma(2-\delta)\lambda_k{a}^\delta}>{0}.
$$

This proves the uniqueness of the solution of the considered problem.

Moreover, one can write an upper bound of $\frac{1}{\Delta}$ by using the last evaluation:

\begin{equation}\label{e33}
	\frac{1}{|\Delta|}\leq\frac{M_1}{(k\pi)^2}+M_2,  \, \, \, \, (M_1, M_2=const).
\end{equation}
Now we find an estimate for $B$ by using Lemma 2.1:

$$
|B|\leq\frac{\lambda_k{p}^{1-\alpha}}{\Gamma(\alpha)}|E_{\delta, 1}(-\lambda_k{a}^\delta)|+\lambda_k{a}^{\delta-1}|E_{\delta, \delta}(-\lambda_k{a}^\delta)|\leq
$$
$$
\leq\frac{\lambda_k{p}^{1-\alpha}}{\Gamma(\alpha)}\frac{M}{1+\lambda_k{a}^\delta}+\lambda_k{a}^{\delta-1}\frac{M}{1+\lambda_k{a}^\delta}\leq
$$
$$
\leq\frac{\lambda_k{p}^{1-\alpha}}{\Gamma(\alpha)}\frac{M}{\lambda_k{a}^\delta}+\lambda_k{a}^{\delta-1}\frac{M}{\lambda_k{a}^\delta}=\frac{Mp^{1-\alpha}}{a^\delta\Gamma(\alpha)}+\frac{M}{a}=\frac{M}{a}\left(1+\frac{p^{1-\alpha}}{a^{\delta-1}}\right)=M_3, \, \, (M_3=const).
$$

Using the last result and (\ref{e33}) we estimate  $\left|\frac{B}{\Delta}\varphi_k\right|$:
$$
\left|\frac{B}{\Delta}\varphi_k\right|\leq\frac{M_1M_3}{(k\pi)^2}|\varphi_k|+\frac{M_2M_3}{k\pi}|\varphi_{1k}|\leq\frac{M_1M_3}{(k\pi)^2}|\varphi_k|+\left(\frac{M_2M_3}{k\pi}\right)^2+|\varphi_{1k}|^2,
$$
here $\varphi_{1k}=2\int_{0}^{1}\varphi'(x)\sin(k\pi{x})dx$.  Now let us find the upper bound of $C$:
$$
|C|\leq\int_{0}^{a}|a-s|^{\delta+q-1}|E_{\delta, \delta+q}(-\lambda_k(a-s)^\delta)||f_k(s)|ds+
$$
$$
+\int_{0}^{a}|a-s|^{\delta+q-2}|E_{\delta, \delta+q-1}(-\lambda_k(a-s)^\delta)||f_k(s)|ds\leq
$$
$$
\leq\int_{0}^{a}|a-s|^{\delta+q-1}\left|\frac{M}{1+\lambda_k|a-s|^\delta}\right||f_k(s)|ds+
$$
$$
+\int_{0}^{a}|a-s|^{\delta+q-2}\left|\frac{M}{1+\lambda_k|a-s|^\delta}\right||f_k(s)|ds\leq\frac{M_4}{(k\pi)^2}, \, \, \, (M_4=const.)
$$
Here we imply that $f_{x}(x,t)\in{L}^1(0, a)$ for convergence of the last integral.

Then, the  estimate for $\frac{C}{\Delta}$ is

$$
\left|\frac{C}{\Delta}\right|\leq\frac{M_1M_4}{(k\pi)^4}+\frac{M_3M_4}{(k\pi)^2}.
$$

Finally, we find the estimate for $|\psi_k|$:
\begin{equation}\label{e34}
|\psi_k|\leq\frac{M_1M_3}{(k\pi)^2}|\varphi_k|+\left(\frac{M_2M_3}{k\pi}\right)^2+|\varphi_{1k}|^2+\frac{M_1M_4}{(k\pi)^4}+\frac{M_3M_4}{(k\pi)^2}<\infty,
\end{equation}
where $\varphi'\in{L}^2(0,1)$.

From (\ref{e32}) and in the same way one can show that $|\tau_k|\leq\frac{M_5}{(k\pi)^2}+|\varphi_{1k}|^2<\infty$,     \, \,  $M_5=const.$

\smallskip

For proving the existence of the solution, we need to show uniform convergency of series representations of  $u(x, t)$, $u_{xx}(x,t)$,  $^C\left(t^\theta\frac{\partial}{\partial{t}}\right)^{\alpha}u(x, t)$ and $D_{t}^{(\gamma, \beta)\mu}u(x, t)$ by using the solution (\ref{e4}) and (\ref{e27}) in $\Omega_1$ and $\Omega_2$ respectively.

In \cite{L22}, the uniform convergence of series of $u(x,t)$ and $u_{xx}(x,t)$ showed for $t>0$. Similarly, for $t>a$, we obtain the following estimate:
$$
|u(x,t)|\leq{M}\sum_{k=1}^{\infty}\Big(\frac{|\tau_k|}{p^\alpha+(k\pi)^2|t^p-a^p|^\alpha}+\frac{1}{(k\pi)^2}\int_{a}^{t}|t^p-\tau^p|^{\alpha-1}f_{2k}(\tau)d(\tau^p)+\\
$$
$$
+\int_{a}^{t}\frac{|t^p-\tau^p|^{2\alpha-1}}{p^\alpha+(k\pi)^2|t^p-a^p|^\alpha}f_{2k}(\tau)d(\tau^p)\Big),
$$
where $f_{2k}(t)=2\int_{0}^{1}f_{xx}(x,t)\sin(k\pi{x})dx$.

Since $|\tau_k|<\infty$ and $f(\cdot, t)\in{C}^3[0,1]$, then the above series converges and hence, by the Weierstrass M-test the series  of $u(x, t)$ is uniformly convergent in $\Omega_1$.

The series of $u_{xx}(x,t)$ is written in the form below
$$
u_{xx}(x,t)=-\sum_{k=1}^{\infty}(k\pi)^2\left(\tau_kE_{\alpha, 1}\left[\frac{(k\pi)^2}{p^\alpha}(t-a)^{p\alpha}\right]+G_k(t)\right)\sin(k\pi{x}).
$$
We obtain the following estimate:
$$
|u_{xx}(x,t)|\leq{M}\sum_{k=1}^{\infty}\Big(\frac{(k\pi)^2|\tau_k|}{p^\alpha+(k\pi)^2|t^p-y^p|^\alpha}+\frac{1}{(k\pi)^2}\int_{a}^{t}|t^p-\tau^p|^{\alpha-1}|f_{4k}(\tau)|d(\tau^p)\\
$$
$$
+\int_{a}^{t}\frac{|t^p-\tau^p|^{2\alpha-1}}{p^\alpha+(k\pi)^2|t^p-\tau^p|^\alpha}|f_{4k}(\tau)|d(\tau^p)\Big),
$$
where $f_{4k}(t)=2\int_{0}^{1}\frac{\partial^4}{\partial{x}^4}f(x,t)\sin(k\pi{x})dx$ and $f(0,t)=f(1,t)=f_{xx}(0,t)=f_{xx}(1,t)=0$.

Since $\tau(0)=\tau(1)=0$ and $\frac{\partial^4{f}}{\partial{x}^4}(\cdot, t)\in{L}^1(0,1)$, then using integration by parts, we arrive at the following estimate

$$
|u_{xx}(x,t)|\leq{M}\sum_{k=1}^{\infty}\left(\frac{1}{k}|\tau_{1k}|+\frac{1}{k^2}\right)\leq\frac{M}{2}\Big(\sum_{k=1}^{\infty}\frac{3}{k^2}+\sum_{k=1}^{\infty}|\tau_{1k}|^2\Big),
$$
where $\tau_{1k}=2\int_{0}^{1}\tau'(x)\sin(k\pi{x})dx$
Then, the Bessel inequality for trigonometric functions implies
$$
|u_{xx}(x,t)|\leq\frac{M}{2}\left(\sum_{k=1}^{\infty}\frac{3}{k^2}+||\tau'(x)||^2_{L^2(0,1)}\right).
$$
Thus, the series in the expression of $u_{xx}(x,t)$ is bounded by a convergent series which is uniformly convergent according to the Weierstrass M-test. Then, the series of $^C(t^\theta\frac{\partial}{\partial{t}})^\alpha{u}(x,t)$ which can be written by
$$
^C\left(t^\theta\frac{\partial}{\partial{t}}\right)^\alpha{u}(x,t)=-\sum_{k=1}^{\infty}(k\pi)^2\left(\tau_kE_{\alpha, 1}\Big[-\frac{(k\pi)^2}{p^\alpha}(t-a)^{p\alpha}\Big]+G_k(t)\right)\sin(k\pi{x})+f(x,t),
$$
has uniform convergence which can be showed in the same way to the uniform convergence of the series of $u_{xx}(x,t)$ (see \cite{L22}).

Now we need to show that the series of $u(x,t)$ and its derivatives should converge uniformly in $\Omega_2$ by using (\ref{e27}). We estimate

$$
|u(x,t)|\leq|\varphi_k||t^{(\beta-2)(1-\mu)}E_{\delta, \delta+\mu(2-\gamma)-1}(-\lambda_k{t}^\delta)|+|\psi_k||t^{\mu+(\beta-1)(1-\mu)}E_{\delta, \delta+\mu(2-\gamma)}(-\lambda_k{t}^\delta)|+
$$
$$
+\int_{0}^{1}|t-\tau|^{\delta-1}|E_{\delta, \delta}\Big(-\lambda_k(t-\tau)^\delta\Big)||f_k(\tau)|d\tau.
$$
Consider estimates of the Mittag-Leffler function (see Lemma 2.1)
$$
|u(x,t)|\leq\frac{|t^{(\beta-2)(1-\mu)}||\varphi_k|M}{1+\lambda_k|{t}^\delta)|}+\frac{|t^{(\beta-1)(1-\mu)}||\psi_k|M}{1+\lambda_k|{t}^\delta)|}
$$
$$
+\int_{0}^{t}|t-\tau|^{\delta-1}\frac{M}{1+\lambda_k|{(t-\tau)}^\delta)|}|f_k(\tau)|d\tau,
$$
where $f_{1k}(t)=\int_{0}^{1}f_{x}(x,t)\sin(k\pi{x})dx$, $f_{x}(\cdot,t)\in{L^1[0,1]}$ and $\varphi'\in{L}^2[0,1]$. Then we obtain the estimate
$$
|u(x,t)|\leq\sum_{k=1}^{\infty}\frac{N_1}{(k\pi)^2}, \, \, \, \, (N_1=\textrm{const}),
$$
for all $t>\bar{t}>0, \, \, \, 0\leq{x}\leq{1}$.

In the similar way one can show that
$$
|u_{xx}(x,t)|\leq\sum_{k=1}^{\infty}(k\pi)^2\Big[|\varphi_k||\frac{M}{1+\lambda_k{t}^\delta}|+(\frac{K_1}{(k\pi)^2}+|\varphi_{1k}|^2)|\frac{M}{1+\lambda_k{t}^\delta}|+
$$
$$
+\int_{0}^{t}|t-\tau|^{\delta-1}\frac{M}{1+\lambda_k{t}^\delta}|k_{2k}(\tau)|d\tau\Big], \, \, \, \, (M=\textrm{const}).
$$
Then, using Bessel's inequality and  $\varphi'\in{L}^2[0,1]$,  $f_{xxx}(\cdot,t)\in{L}^1(0,1)$, we get
$$
|u_{xx}(x,t)|\leq\frac{1}{2}\Big(\sum_{k=1}^{\infty}\frac{N_2}{(k\pi)^2}+||\varphi'(x)||^2_{L^2[0,1]}\Big).
$$
We have also used $2ab\leq{a^2+b^2}.$

Using the equation in $\Omega_2$, we write $D_{t}^{(\gamma, \beta)\mu}u(x,t)$ in the form
$$
D_{t}^{(\gamma, \beta)\mu}u(x,t)=u_{xx}(x,t)+f(x,t)$$
and its uniform convergence can be done in a similar way to the uniform convergence of $u_{xx}(x,t)$ as
$$
|D^{(\gamma, \beta)\mu}_tu(x,t)|\leq\sum_{k=1}^{\infty}\frac{N_3}{(k\pi)^2}, \, \, \, \, (N_3=\textrm{const}).
$$

Finally, considering the Weierstrass M-test, the above arguments prove that Fourier series in (\ref{e4}) and (\ref{e27})  converge uniformly in the domains $\Omega_1$ and $\Omega_2$. This is the proof that the considered problem's solution exists in $\Omega$.  Theorem 3.1 is proved.

\section*{Appendix}

Here we write derivation of the series $^{C}\left(t^{\theta}\frac{\partial}{\partial{t}}\right)^\alpha{u}(x,t)$ in (\ref{e4}). Using  relation (\ref{eq3}) we get:
$$
^{C}\left(t^{\theta}\frac{\partial}{\partial{t}}\right)^\alpha{u}(x,t)=\sum_{k=0}^{\infty}\Big[(t^\theta\frac{\partial}{\partial{t}})^\alpha\left(\tau_kE_{\alpha,1}\Big[-\frac{(k\pi)^2}{p^\alpha}(t^p-a^p)^\alpha\Big]+G_k(t)\right)
-\frac{\tau_k(t^p-a^p)^\alpha}{p^{-\alpha}\Gamma(1-\alpha)}\Big]\sin(k\pi{x}).
$$
The hyper-Bessel derivative of the Mittag-Leffler function is
$$
\Big(t^\theta\frac{\partial}{\partial{t}}\Big)^\alpha\tau_kE_{\alpha,1}\left(-\frac{(k\pi)^2}{p^\alpha}(t^p-a^p)^\alpha\right)=\tau_kp^\alpha(t^p-a^p)^{-\alpha}E_{\alpha, 1-\alpha}\big[-\lambda(t^p-a^p)^\alpha\big].
$$
Using the Lemma 2.3, we can write the last expression as follows
$$
\Big(t^\theta\frac{\partial}{\partial{t}}\Big)^\alpha\tau_kE_{\alpha,1}\Big[-\frac{(k\pi)^2}{p^\alpha}(t^p-a^p)^\alpha\Big]=\frac{\tau_kp^\alpha(t^p-a^p)^{-\alpha}}{\Gamma(1-\alpha)}+\tau_k(k\pi)^2E_{\alpha, 1}\big[-\frac{(k\pi)^2}{p^\alpha}(t^p-a^p)^\alpha\big].
$$

Then evaluating $\Big(t^\theta\frac{\partial}{\partial{t}}\Big)^\alpha{G}_k(t)$ gives that
$$
\Big(t^\theta\frac{\partial}{\partial{t}}\Big)^\alpha{G}_k(t)=\Big(t^\theta\frac{\partial}{\partial{t}}\Big)^\alpha\Big(f^*_k(t)+\lambda^*\int_{a}^{t}(t^p-a^p)^{\alpha-1}E_{\alpha, \alpha}\big[\lambda^*(t^p-a^p)\big]f^*_k(\tau)d(\tau^p)\Big)=
$$
$$
=p^\alpha{t}^{-p\alpha}D_{p, a+}^{-\alpha, \alpha}\Big(\frac{1}{p^\alpha}I_{p, a+}^{-\alpha, \alpha}t^p\alpha{f}_k(t)+\lambda^*\int_{a}^{t}(t^p-a^p)^{\alpha-1}E_{\alpha, \alpha}\big[\lambda^*(t^p-a^p)\big]f^*_k(\tau)d(\tau^p)\Big)=
$$
$$
=f_{k}(t)+p^{\alpha}t^{-p\alpha}D_{p, a+}^{-\alpha, \alpha}\Big(\lambda^*\int_{a}^{t}(t^p-a^p)^{\alpha-1}E_{\alpha, \alpha}\big[\lambda^*(t^p-a^p)\big]f^*_k(\tau)d(\tau^p)\Big),
$$
where $\lambda^*=-\frac{\lambda_k}{p^\alpha}$ and $f^*_k(t)=\frac{1}{p^\alpha\Gamma(\alpha)}\int_{a}^{t}(t^p-\tau^p)^{\alpha-1}f_k(\tau)d(\tau^p)$.

The second term in the last expression can be simplified using the Erd'elyi-Kober fractional derivative for $n=1$,
$$
-\lambda_kt^{-p\alpha}\left(1-\alpha+\frac{t}{p}\frac{d}{dt}\right)\frac{t^{-p(1-\alpha)}}{\Gamma(1-\alpha)}\int_{a}^{t}(t^p-\tau^p)^{-\alpha}d(\tau^p)\int_{a}^{\tau}(\tau^p-s^p)^{\alpha-1}E_{\alpha, \alpha}\big[\lambda^*(\tau^p-s^p)^\alpha\big]f^*_k(s)d(s^p)=
$$
$$
-\lambda_kt^{-p\alpha}\left(1-\alpha+\frac{t}{p}\frac{d}{dt}\right)\frac{t^{-p(1-\alpha)}}{\Gamma(1-\alpha)}\int_{a}^{t}f^*_k(s)d(s^p)\int_{s}^{t}(t^p-\tau^p)^{-\alpha}(\tau^p-s^p)^{\alpha-1}E_{\alpha, \alpha}\big[\lambda^*(\tau^p-s^p)^\alpha\big]d(\tau^p)=
$$
$$
-\lambda_kt^{-p\alpha}\left(1-\alpha+\frac{t}{p}\frac{d}{dt}\right)t^{-p(1-\alpha)}\int_{a}^{t}E_{\alpha, 1}\left[\lambda^*(t^p-s^p)^{\alpha}\right]f^*_k(s)d(s^p)=
$$
$$
-\lambda_k(1-\alpha)t^{-p}\int_{a}^{t}E_{\alpha, 1}\left[\lambda^*(t^p-s^p)^{\alpha}\right]f^*_k(s)d(s^p)-
$$
$$
-\frac{\lambda_kt^{-p\alpha+1}}{p}\frac{d}{dt}\left(t^{-p(1-\alpha)}\int_{a}^{t}E_{\alpha, 1}\left[\lambda^*(t^p-s^p)^{\alpha}\right]f^*_k(s)d(s^p)\right)=
$$
$$
=-\lambda_k(1-\alpha)t^{-p}\int_{a}^{t}E_{\alpha, 1}\left[\lambda^*(t^p-s^p)^{\alpha}\right]f^*_k(s)d(s^p)+
$$
$$
+\lambda_k(1-\alpha)t^{-p}\int_{a}^{t}E_{\alpha, 1}\left[\lambda^*(t^p-\tau^p)^{\alpha}\right]f^*_k(\tau)d(\tau^p)-\lambda_kf^*_k(t)-
$$
$$
-\lambda_kt^{1-p}\int_{a}^{t}\lambda^*(t^p-\tau^p)^{\alpha-1}E_{\alpha, \alpha}\left[\lambda^*(t^p-\tau^p)^{\alpha}\right]f^*_k(\tau)d(\tau^p)=
$$
$$
=-\lambda_k\left(f^*_k(t)+\lambda^*\int_{a}^{t}(t^p-a^p)^{\alpha-1}E_{\alpha, \alpha}\big[\lambda^*(t^p-a^p)\big]f^*_k(\tau)d(\tau^p)\right)=-\lambda_kG_k(t).
$$
Hence, we get
$$
^{C}\left(t^{\theta}\frac{\partial}{\partial{t}}\right)^\alpha{u}(x,t)=-\sum_{k=0}^{\infty}(k\pi)^2\left[\tau_kE_{\alpha, 1}\left(\frac{(k\pi)^2}{p^\alpha}(t^p-a^p)\right)+G_k(t)\right]\sin(k\pi{x})+f(x,t).
$$
This proves that solution (\ref{e4}) satisfies the equation
$$
^C\Big(t^\theta\frac{\partial}{\partial{t}}\Big)^{\alpha}u(x,t)-u_{xx}(x,t)=f(x,t).
$$

We would like to note that using the result of this work, one can consider FPDE with the Bessel operator considering local \cite{krma} and non-local boundary value problems \cite{karruzh}. In that case the Fourier-Bessel series will play an important role. The other possible applications are related with the consideration of more general operators in space variables. For instance, in \cite{rtb1}, very general positive operators have been considered, and the results of this paper can be extended to that setting as well.

\section{Acknowledgement}

The second author was partially supported by the FWO Odysseus 1 grant G.0H94.18N: Analysis and Partial Differential Equations and by the Methusalem programme of the Ghent University Special Research Fund (BOF) (Grant number 01M01021).

\end{document}